\documentclass[12pt,reqno]{amsart}
\usepackage{enumerate, latexsym, amsmath, amsfonts, amssymb, amsthm, color}
\def\pmod #1{\ ({\rm{mod}}\ #1)}
\def\Z{\Bbb Z}
\def\N{\Bbb N}

\def\l{\left}
\def\r{\right}
\def\bg{\bigg}
\def\({\bg(}
\def\){\bg)}
\def\t{\text}
\def\f{\frac}

\def\ls{\leqslant}
\def\gs{\geqslant}
\def\se {\subseteq}
\def\sm{\setminus}
\def\bi{\binom}

\def\eq{\equiv}

\def\Proof{\noindent{\it Proof}}

\theoremstyle{plain}
\newtheorem{theorem}{Theorem}

\newtheorem{lemma}{Lemma}

\theoremstyle{definition}

\theoremstyle{remark}

 \vspace{4mm}

\begin{document}

\hbox{Preprint}

\title
[On an inverse problem for restricted sumsets]
{On an inverse problem \\ for restricted sumsets}

\author
[] {Xin-Qi Luo and Zhi-Wei Sun}

\address {(Xin-Qi Luo) Department of Mathematics, Nanjing
University, Nanjing 210093, People's Republic of China}
\email{lxq15995710087@163.com}

\address{(Zhi-Wei Sun, corresponding author) Department of Mathematics, Nanjing
University, Nanjing 210093, People's Republic of China}
\email{zwsun@nju.edu.cn}

\keywords{Additive combinatorics, restricted sumset, inverse problem.
\newline \indent This research is supported by the National Natural Science Foundation of China (grant 11971222).}

\begin{abstract}
Let $n$ be a positive integer, and let $A$ be a set of $k\gs2n-1$ integers. For the restricted sumset
$$
S_n(A)=\{a_1+\cdots +a_n:\ a_1,\ldots,a_n\in A,\ \t{and}\  a_i^2\neq a_j^2\ \t{for} \ 1\ls i<j\ls n\},
$$
by a 2002 result of Liu and Sun we have
$$|S_n(A)|\gs (k-1)n-\f 32n(n-1)+1.$$
In this paper, we determine the structure of $A$ when the lower bound is attained.
\end{abstract}
\maketitle

\section{Introduction}
\setcounter{lemma}{0}
\setcounter{theorem}{0}
\setcounter{corollary}{0}
\setcounter{remark}{0}
\setcounter{equation}{0}

Let $G$ be an additive abelian group. For any finite $A\se G$ and $n\in\Z^+=\{1,2,3,\ldots\}$,
we define
$$n^\wedge A=\{a_1+a_2+\cdots +a_n:\ a_1,\ldots,a_n\in A, \ \t{and}\ a_i\not=a_j\ \ \t{for} 1\ls i<j\ls n\}.$$
When $G$ is the additive group $\Z$ of integers, it is known that
\begin{equation}\label{nA} |n^\wedge A|\gs n|A|-n^2+1=n(|A|-n)+1
\end{equation}
(cf. \cite{N}).

Let $A$ be a set of $k$ integers. The inverse problem for $n^\wedge A$ studies
the structure of $A$ when equality in \eqref{nA} holds.
M.B. Nathanson \cite{N95} noted that if $k=4$ and $A=\{a_0,a_1,a_2,a_3\}$
with $a_3-a_2=a_1-a_0$ then the equality $|2^\wedge A|= 2k-2^2+1$ holds but $A$ may not be an arithmetic progression (AP).
He also obtained the following result.

\begin{theorem}[Nathanson \cite{N95}]\label{Nath} Let $A$ be a set of $k\gs 5$ integers, and let $n\in\{2,\ldots,k-2\}$. Then
$|n^\wedge A|= nk-n^2+1$ if and only if $A$ is an AP.
\end{theorem}

Let $p$ be a prime, and let $A\se\Z/p\Z$.
In 1964 P. Erd\H os and H. Heilbronn \cite{EH64} conjectured the inequality
 $|2^\wedge A|\gs \min \{p, 2|A|-3\}$.
In 1994, Dias da Silva and H. Hamidoune \cite{DH94} confirmed this conjecture
via exterior algebra, and they proved further that
 $$|n^\wedge A|\gs\min \{p,n|A|-n^2+1\}$$
 for any positive integer $n$. A simple proof of this via the so-called polynomial method
 was given by N. Alon, Nathanson and I.Z. Ruzsa \cite{ANR1,ANR2}.
 In 2005, G. K\'arolyi \cite{K05} used the Combinatorial Nullstellensatz of Alon \cite{A99}
 to prove that if $|A|=k\gs 5$, $p>2k-3$ and $|2^\wedge A|=2k-3$ then $A$ is an AP, which was also deduced by
 K\'arolyi and R. Paulin \cite{KP13} in 2013 via exterior algebra.

Let $A$ be a $k-$element subset of a field $F$ of characteristic $p$ as $\{a_1,a_2\dots a_k\}$. Then $|2^\wedge A|\gs \min \{p,2k-3\}$.
We define
$$S_n(A)=\{a_1+\dots +a_n:a_i\in A_i,\ \t{and} \ a_i^2\neq a_j^2\ \t{if} \ i\neq j\}.$$
If we add condition $a_i^2\neq a_j^2$, then we only should remove element $0=a_i+a_j\in 2^\wedge A$ by $a_i=-a_j$ for $1\ls i,j \ls k$, hence we have $|S_2(A) |\gs \min \{p-1,2k-4\}$.

In 2002, J.-X. Liu and Z.-W. Sun \cite{LS02} used the polynomial method to deduce the following result
which extends the da Silva--Hamidoune theorem.

\begin{theorem} [Liu and Sun \cite{LS02}] \label{Th1.1} Let $k,m,n$ be positive integers with $k>m(n-1)$, and let $F$ be a field of characteristic $p$ where $p$ is zero or greater than $K=(k-1)n-(m+1)\binom n2$. Let $A_1, A_2\dots A_n$ be subsets of $F$ for which
$$|A_n|=k \ and\ |A_{i+1}|-|A_i|\in \{0,1\} \ for \ i=1,\dots ,n-1,$$
If $P(x)\in F[x]$ is monic and of degree $m$, then we have
$$|\{a_1+\dots +a_n:a_i\in A_i,\ and\ P(a_i)\neq P(a_j)\ if\ i\neq j\}|\gs K+1.$$
\end{theorem}

As a consequence of this result, if the characteristic of a  field $F$ is zero or greater than
 $(k-1)n-3\bi n2$ with $k\gs 2n-1$, then for any finite $A\se F$ with $|A|=k$ we have
\begin{equation}\label{Sn} |S_n(A)|\gs (k-1)n-3\bi n2+1,
\end{equation}
where
\begin{equation} S_n(A):=\{a_1+\cdots+a_n:\ a_i\in A,\ \t{and}\ a_i^2\not=a_j^2
\ \t{if}\ i\not=j\}.
\end{equation}
In 2012 \'E. Balandraud \cite{Ba12} applied this result to prove that
 if $p$ is a prime and $A\se \Z/p\Z$ with $a+b\not=0$ for all $a,b\in A$ then
$$\bg|\bg\{\sum_{a\in B}a:\ \emptyset\not=B\se A\bg\}\bg|\gs\min\l\{p,\,\f{|A|(|A|+1)}2\r\}.$$
This implies that for any prime $p$ we have
\begin{align*}&\max\bg\{|A|:\ \sum_{a\in B}a\not=0\ \t{for any}\ \emptyset\not=B\se A\bg\}
\\=&\max\bg\{k\in\Z:\ \f{k(k+1)}2<p\bg\}=\l\lfloor\f{\sqrt{8p-7}-1}2\r\rfloor,
\end{align*}
which is a conjecture of Erd\H os and J.L. Selfridge (cf. \cite[p.\,194]{G94}).

Now, let $A\se\Z$ with $|A|=k\in\Z^+$. Then we have the inequality \eqref{Sn} if $k\gs 2n-1$.
Note that the inequality \eqref{Sn} may fail if $k<2n-1$, for example, $S_n(A)=\emptyset$
when $$
|A|\se\{\pm1,\pm2,\ldots,\pm(n-1)\}.
$$
In this paper we study the inverse problem for $S_n(A)$. Namely, when $k\gs 2n-1$ we determine completely when the lower bound in \eqref{Sn} is attained.

Now we state our main results.

\begin{theorem}\label{Th1.2}
Let $A$ be a set of finite integers with $|A|\gs 3$. Then
$$|S_2(A)|=(|A|-1)2-3\bi 22+1=2|A|-4,$$
if and only if one of the following {\rm (i)-(v)} holds.

{\rm (i)} $A=\{c,\pm d\}$ for some $c\in\Z$ and $d\in\Z^+$ with $|c|\neq d$.

{\rm (ii)} $A=\{\pm c,\pm d\}$ for some $c,d\in\Z^+$ with $c\neq d$.

{\rm (iii)}$A=\{\pm d\}\cup \{c\pm d\}$ for some $c\in\Z$ and $d\in\Z^+$ with $c\neq 0,\pm 2d$.

{\rm (iv)} $A=\{rd:\ s\ls r\ls t\}$ for some integers $d\not=0$, $s\ls-1$ and $t\gs1$ with $t-s\gs4$.

{\rm (v)} $A=\{(2r-1)d:\ s\ls r\ls t\}$ for some integers $d\not=0$, $s\ls0$ and $t\gs 1$
with $t-s\gs4$.
\end{theorem}

\begin{theorem}\label{Th1.3} Let $n>2$ be an integer, and let $A\se\Z$ with $|A|=k\in\{2n-1,2n\}$.
Concerning the equality
\begin{equation}\label{A=}|S_n(A)|= (k-1)n-3\binom n2+1=kn-\frac 32 n^2+\frac n2+1,
\end{equation}
we have the following results.

{\rm (i)} Assume that $k=2n-1$. Then \eqref{A=} holds, if and only if there are $b,d\in\Z^+$ and $c\in\Z$
such that $A=\{c,\pm d\}\cup\{\pm (b+jd):\ j=0,\ldots,n-2\}$, and $b=d$ if $n>3$.

{\rm (ii)} Suppose that $k=2n$.  Then \eqref{A=} holds, if and only if
there are $b,d\in\Z^+$ such that $A=\{\pm d\}\cup\{\pm (b+jd):\ j=0,\ldots,n-1\}$, and
$b=d$ if $n>3$.
\end{theorem}

{\it Example} 1.1. For $A=\{\pm1,\pm4,6\},\{\pm1,\pm3,\pm4\}$, we have
$$|S_3(A)|=(|A|-1)3-3\bi 32+1=3|A|-11;$$
in fact,
$S_3(\{\pm1,\pm4,6\}) =\{1,3,9,11\}$ and
$$S_3(\{\pm1,\pm3,\pm4\}) =\{0,\pm 2,\pm 6,\pm 8\}.$$
For any integer $n>3$, if $A=\{\pm 1,\pm 2,\ldots, \pm n\}$ then
$$S_n(A) =\l\{\f{n(n+1)}2-2j:\ j=0,\ldots,\f{n(n+1)}2\r\}$$
and hence
$$|S_n(A)|=(|A|-1)n-3\bi n2+1.$$

\begin{theorem}\label{Th1.4} Let $n>2$ and $A$ be a set of $k$ integers with $k\gs 2n+1$. Then
$$|S_n(A)|> (k-1)n-3\binom n2+1.$$
\end{theorem}

{\it Example} 1.2. For $A=\{\pm 1,\pm 3,\pm 5,7\},$ we have
$$
S_3(A)=\{\pm 1,\pm 3,\pm 7,\pm 9,5 ,11 ,13,15 \}
$$
and hence $|S_3(A)|=12>(|A|-1)3-3\binom 32+1=10$

We are going to prove Theorem 1.3 and Theorem 1.4(i) in Sections 2 and 3, respectively.
Theorem 1.4(ii) and Theorem 1.5 will be proved in Sections 4.

\section{Proof of Theorem 1.3}
\setcounter{lemma}{0}
\setcounter{theorem}{0}
\setcounter{corollary}{0}
\setcounter{remark}{0}
\setcounter{equation}{0}

\noindent{\bf Proof of Theorem 1.3.} The ``if" direction is easy since we can verify the result directly.

Now we assume that $|S_2(A)|=2|A|-4$. As $|2^{\wedge} A|\gs 2|A|-3$ by \eqref{nA}, we have $S_2(A)\not=2^{\wedge} A$
and hence $-A\cap(A\sm\{0\})\neq \emptyset$. It follows that
\begin{equation} \label{-A} -A\cap(A\sm\{0\})=\{\pm d:\ d\in D\}\ \ \t{for some}\ D\se\Z^+.
\end{equation}

{\it Case} 1. $|A|=3.$

In view of \eqref{-A}, for some $d\in\Z^+$ we have $\pm d\in A$. Thus $A=\{c,\pm d\}$
for some $c\in\Z$ with $c\not=\pm d$.
\medskip

{\it Case} 2. $|A|=4.$

When $|-A\cap(A\sm\{0\})|=4$, by \eqref{-A} we have $A=\{\pm c,\pm d\}$ for some $c,d\in\Z^+$ with $c\not=d$.

Now suppose that $-A\cap(A\sm\{0\})=\{\pm d\}$ for some $d\in\Z^+$.
Write $A=\{\pm d\}\cup \{c-d,b\}$ with $b,c\in\Z$ and $c-d<b$. As $|A|=4$, we see that $S_2(A)$ contains
the five numbers
$$(c-d)+(-d),\ b+(-d),\ (c-d)+d,\ b+d,\ b+(c-d).$$
Note that $b+c-d$ is different from all the four numbers $c-2d,b-d,c,b+d$. Also,
$$c-2d<b-d<b+d\ \ \t{and}\ \ c-2d<c<b+d.$$
Since $|A|=4$, we must have $b-d=c$. Thus $A=\{\pm d,c\pm d\}.$
\medskip

{\it Case} 3. $|A|\gs 5.$

If $0\not\in 2^\wedge A$, then all those $a^2\ (a\in A)$ are pairwise distinct, and hence $S_2(A)=2^\wedge A$.
 With the aid of \eqref{nA}, we have $|2^{\wedge} A|>2|A|-4=|S_2(A)|$ and hence $2^\wedge A\sm S_2(A)$ contains $0$.
 If $c$ is a nonzero element of $2^\wedge A$, then $c\in S_2(A)$ since $c=a+b$ for some $a,b\in A$ with $a\pm b\not=0$.
 Therefore $2^\wedge A=S_2(A)\cup\{0\}$ and
$|2^{\wedge} A|= 2|A|-3$. Applying Theorem 1.1, we see that $A$ is an AP.
Let $d$ be the least positive integer in $A$, and let $q\in\Z^+$ be the common difference of $A$.

If $0\in A$, then $q=d-0=d$ and hence
$$A=\{rd:\ s\ls r\ls t\}$$
for some integers $d\not=0$, $s\ls-1$ and $t\gs1$.

Now we assume $0\notin A$.
In view of \eqref{-A},  we have $\pm a\in A$ for some $a\in A$, and hence $d\eq-d\pmod q$.
By the choice of $d$ and $q$, we have $d-q\ls0$.
As $0\not\in A$, we have  $d<q\ls 2d$ and hence $q=2d$ since $q\mid 2d$.
Thus
$$A=\{(2r-1)d:\ s\ls r\ls t\}$$
for some integers $d\not=0$, $s\ls0$ and $t\gs 1$.

Combining the above, we have completed the proof of Theorem 1.3. \qed

\section{Proof of Theorem 1.4(i)}
\setcounter{lemma}{0}
\setcounter{theorem}{0}
\setcounter{corollary}{0}
\setcounter{remark}{0}
\setcounter{equation}{0}

For the sake of convenience, for a finite nonempty subset $S$ of $\Z$, we define
$$\sigma(S)=\sum_{s\in S}s,\ \min S=\min\{s:\ s\in S\}
\ \t{and}\  \max S=\max\{s:\ s\in S\}.$$

\begin{lemma}\label{lemma1}
Let $n>2$ and $k\gs 2n-1$ be integers. Suppose that $A$ is a set of $k$ integers satisfying \eqref{A=}.
Then $A$ contains at least $n-1$ negative integers and also at least $n-1$ positive integers.
\end{lemma}
\Proof. If $A\se\N$, then
$$|S_n(A)|=|n^\wedge A|\gs kn-n^2+1>kn-\f 32n^2+\f n2+1$$
by \eqref{nA}. As \eqref{A=} holds, we must have $A\not\se \N$.
Similarly, $-A\not\se\N$ since $|S_n(-A)|=|S_n(A)|$. As
$$|\{a\in A:\ a<0\}|+|\{a\in A:\ a>0\}|\gs k-1\gs 2n-2,$$
without loss of generality we may assume that $|\{a\in A:\ a<0\}|\ls n-1$.
Write all the $k$ elements of $A$ in the ascending order
$$b_r<b_{r-1}<\cdots <b_1<a_1<\cdots <a_{k-r}$$
with $b_1<0$ and $r\ls n-1$.
For $\emptyset\not=S\se A$ with  $-S\cap (S\sm\{0\})=\emptyset$, when $b_i\notin S$ we define
$$
U_1(S,b_i)=\{b_i\}\cup (S\sm\{u(S,b_i)\}),
$$
where
$$u(S,b_i)=\begin{cases}\min S&\t{if}\ \min S>-b_i,
\\ \max\{s\in S:\ s\ls -b_i\}&\t{otherwise}.
\end{cases}$$
For $2\ls m\ls n-1$ and $m$ distinct elements $b_1',\ldots,b_m'$ of $A\sm S$, we define
$$
U_m(S,b'_1,\ldots ,b'_m)=\{b'_m\}\cup(U_{m-1}(S,b'_1,\ldots b'_{m-1})\sm\{u(X_m,b_m')\})
$$
where $$X_m=U_{m-1}(S,b'_1,\ldots b'_{m-1})\sm \{b'_1\ldots b'_{m-1}\}.$$

As $k\gs2n-1$ and $r\ls n-1$, we have $k-r\gs n$.

Set $T=\{a_1,\ldots,a_n\}$, and assume that there is an element $c$ of $T$ with $c=0$ or $c\notin -A.$
When $-T\cap (T\sm\{0\})=\emptyset$, it is easy to see that $S_n(A)$ contains the following numbers:
\begin{equation}\label{list}\begin{aligned}&\sigma(T),c+\sigma(U_1(T\sm \{c\},b_i))\ (i=1,\ldots,r),
\\&c+\sigma(U_2(T\sm \{c\},b_r,b_j))\ (j=1,\ldots,r-1),
\\&\cdots,c+\sigma(U_{r}(T\sm \{c\},b_r,b_{r-1}\ldots ,b_1)).
\end{aligned}
\end{equation}
Clearly, $\sigma(T)>c+\sigma( U_1(T\sm \{c\},b_1)).$

Let $1\ls m\ls r-1$ and $1\ls i\ls r-m$. Then
\begin{align*}&\sigma(U_{m+1}(T\sm \{c\},b_r\ldots b_{r-m+1},b_i))
\\=\ &\sigma(U_{m}(T\sm \{c\},b_r\ldots b_{r-m+1}))+b_i-u(X_m,b_{i}),
\end{align*}
where $X_m=U_{m}(T\sm \{c\},b_{r},\ldots b_{r-m})\sm \{b_{r}\ldots b_{r-m}\}.$
Note that $b_i<\min X_m$ since $u(X_m,b_i)\in X_m$.
Thus
$$\sigma(U_{m}(T\sm \{c\},b_r\ldots b_{r-m+1}))>\sigma(U_{m+1}(T\sm \{c\},b_r\ldots b_{r-m+1},b_i)).$$

Next, we show that
\begin{align*}
&\sigma(U_m(T\sm \{c\},b_r\ldots b_{r-m+2},b_i))
>\sigma(U_m(T\sm \{c\},b_r\ldots b_{r-m+2},b_{i+1})).
\end{align*}
Observe that
$$\sigma(U_m(T\sm \{c\},b_r\ldots b_{r-m+2},b_i))
-\sigma(U_m(T\sm \{c\},b_r\ldots b_{r-m+2},b_{i+1}))$$
equals $b_i-b_{i+1}+u(X_m,b_{i+1})-u(X_m,b_{i})$, which is clearly greater than $u(X_m,b_{i+1})-u(X_m,b_{i})$.
So we it suffices to prove that $u(X_m,b_{i+1})\gs u(X_m,b_{i})$.

{\it Case} 1. $\min X_m>-b_{i+1}$.

In this case, $\min X_m>-b_{i}$ and hence
$$u(X_m,b_{i+1})=\min X_m=u(X_m,b_{i}).$$

{\it Case} 2. $\min X_m\ls -b_{i+1}$.

In this case,
$$u(X_m,b_{i+1})=\max\{s\in X_m:\ s\ls -b_{i+1}\}.$$
If
$$u(X_m,b_{i})=\max\{s\in X_m:\ s\ls -b_i\},$$
then
\begin{align*}u(X_m,b_{i+1})&=\max\{s\in X_m:\ s\ls -b_{i+1}\}
\\&\gs \max\{s\in X_m:\ s\ls -b_i\}=u(X_m,b_i).
\end{align*}
When $u(X_m,b_{i})=\min X_m,$ we have
$$u(X_m,b_{i+1})=\max\{s\in X_m:\ s\ls -b_{i+1}\}\gs \min X_m=u(X_m,b_i).$$

By the above, the $r(r+1)/2+1$ numbers in \eqref{list} are in decreasing order.

Note that we may simply assume $a_1\gs0$.
 Suppose that $r<n-1$. Then $a_1^2,\ldots,a_{k-r}^2$ are pairwise distinct, and hence $\sigma(T)\in S_n(A)$.
Also, $c=0$ or $-c\notin A$ for some $c\in T$, and hence the $r(r+1)/2+1$ numbers in \eqref{list}
are distinct elements of $S_n(\{b_r,\ldots b_1,a_1\ldots a_n\})$.
Note that $S_n(\{a_1,a_2\ldots a_{k-r}\})\gs n(k-r-n)+1$ by \eqref{nA}. Also,
$$S_n(A)\supseteq (S_n(\{b_r,\ldots b_1,a_1\ldots a_n\})\cup S_n(\{a_1,a_2\ldots a_{k-r}\}))$$
and
$$S_n(\{b_r,\ldots b_1,a_1\ldots a_n\})\cap S_n(\{a_1,a_2\ldots a_{k-r}\})=\{a_1+\cdots+a_n\}=\{\sigma(T)\}.$$
It follows that
$$|S_n(A)|\gs\f{r(r+1)}2+1+(n(k-r-n)+1)-1=\f{r(r+1)}2-rn+kn-n^2+1.$$
Let $f_n(s)=s(s+1)/2-sn$. Then
$$f_n(n-2)-f_n(r)=\sum_{r\ls s<n-2}(f_n(s+1)-f_n(s))=\sum_{r\ls s<n-2}(s+1-n)\ls 0.$$
Therefore,
$$|S_n(A)|\gs f_n(n-2)+kn-n^2+1=kn-\f32n^2+\f n2+2$$
which contradicts with the condition $|S_n(A)|=kn-\f 32n^2+\f n2+1$.

By the above, $A$ contains at least $n-1$ negative numbers. Similarly, $-A$ also contains at least $n-1$ negative numbers, i.e., $|A\cap\Z^+|\gs n-1$. This ends our proof. \qed

\medskip
\noindent{\it Proof of Theorem} 1.4(i).
The ``if" direction is easy since we can verify the result directly.

Now assume the equality \eqref{A=} with $k=|A|=2n-1$. Note that $|S_n(-A)|=|S_n(A)|$.
By Lemma 3.1, either $A$ or $-A$ contains exactly $n-1$ negative integers
and $n$ nonnegative integers. Without loss of generality, we simply assume that $A$
consists of $n-1$ negative integers $b_{n-1}<\ldots<b_1$  and $n$ nonnegative integers.
Clearly, for some $c\in A$ with $c\gs0$ we have $c=0$ or $-c\not\in A$.
Let $a_1<\ldots<a_{n-1}$ be the list of elements of $A\sm\{c\}$ in ascending order.

Let $T=\{a_1,\ldots ,a_{n-1},c\}$. By the reasoning concerning \eqref{list} in the proof of Lemma \ref{lemma1}, the
$n(n-1)/2+1$ numbers
\begin{equation}\label{list'}\begin{aligned}&\sigma(T),c+\sigma(U_1(T\sm \{c\},b_i))\ (i=1,\ldots,n-1),
\\&c+\sigma(U_2(T\sm \{c\},b_{n-1},b_j))\ (j=1,\ldots,n-2),
\\&\cdots,c+\sigma(U_{n-1}(T\sm \{c\},b_{n-1},b_{n-2}\ldots ,b_1)).
\end{aligned}\end{equation}
are distinct elements of $S_n(A)$. As \eqref{A=} holds, this list gives all the elements of $S_n(A)$.

We claim that $-b_j\in A$ for all $j=1,\ldots,n-1$. Assume on the contrary that if there a $j\in\{1,\ldots,n-1\}$ such that $-b_j\notin A.$
When $c<a_{n-1},$
for the set $T_1=\{b_j,a_1,\ldots a_{n-2},c\}$,the $(n-1)(n-2)/2+1$ numbers
\begin{align*}&\sigma(T_1),b_j+\sigma(U_1(T_1\sm \{b_j\},b_i))\ (i=1,\ldots ,j-1,j+1,\ldots ,n-1),
\\&b_j+\sigma(U_2(T_1\sm \{b_j\},b_{n-1},b_i))\ (i=2,\ldots ,j-1,j+1,\ldots ,n-1),
\\&\cdots,b_1+\sigma(U_{n-2}(T_1\sm \{b_j\},b_{n-1},\ldots ,b_{j+1},b_{j-1},\ldots,b_1)).
\end{align*}
are distinct elements of $S_n(A)$ and they are in descending order,
also the $n$ numbers
$$b_j+a_1+\cdots +a_{n-1}, \cdots \ , a_1+a_2+\cdots +a_{n-1}+c$$
are distinct elements of $S_n(A)$ greater than $\sigma(T_1)$.
As $$\frac {(n-1)(n-2)}2+n+1>\frac {n(n-1)}2+1=|S_n(A)|,$$
we must have $c>a_{n-1}$. Similarly, using
$$T_2=(T_1\sm\{c\})\cup\{a_{n-1}\}=\{b_j,a_1,\ldots a_{n-2},a_{n-1}\}$$ instead of $T_1$
in the above arguments, we also deduce a contradiction.

By the last paragraph, $-b_1<\ldots<-b_{n-1}$ are positive elements of $A$ which are different from $c$ since $-c\not\in A$. Thus $a_i=-b_i$ for $1\ls i\ls n-1$.

We now show that  $a_{i+1}-a_i=a_1$.
Recall that all the elements of $S_n(A)$ are exhibited in the list \eqref{list'}, and
\begin{equation}\label{n1}
c+\sigma(U_1(T\sm \{c\},b_i))\ (i=1,\ldots,n-1),c+\sigma( U_2(T\sm \{c\},b_{n-1},b_1))
\end{equation}
are $n$ consecutive terms on the list \eqref{list'} in descending order. Note also that the $n$ numbers
\begin{equation}\label{n2}
c+\sigma(U_1(T\sm \{c\},b_i))\ (i=1,2),c+\sigma( U_2(T\sm \{c\},b_{i},b_1))\ (i=2,\ldots,n-1 )
\end{equation}
in descending order also belong to $S_n(A)$. As \eqref{n1} and \eqref{n2} have the same first term
and the same last term, the two lists \eqref{n1} and \eqref{n2} are identical. Therefore, for any $2\ls i\ls n-2$, we have
$$c+\sigma(U_1(T\sm \{c\},b_{i+1}))=c+\sigma( U_2(T\sm \{c\},b_{i},b_1))$$
and hence $a_{i+1}-a_i=a_1$.

Note that $d=a_1$ and $b=a_2-a_1$ are positive integers. For each $i=2,\ldots,n-1$ we have $-b_i=a_i=b+(i-1)d$.
Thus
$$A=\{c,\pm d\}\cup\{\pm (b+jd):\ j=1,\ldots,n-2\}.$$

Now assume  $n\gs4$. We want to prove that $b=d$. Note that
\begin{align*}
&c+\sigma( U_1(T\sm \{c\},-a_{n-1}))\\>&\ c+\sigma(U_2(T\sm \{c\},-a_{n-2},-a_2))\\
=&\ d-(b+d)+\sum_{1<j<n-3}(b+jd)-(b+(n-3)d)+b+(n-2)d+c
\end{align*}
and also
$$c+\sigma( U_2(T\sm \{c\},-a_{n-2},-a_2))>c+\sigma( U_2(T\sm \{c\},-a_{n-1},-a_{n-2})).$$
Clearly, $c+\sigma( U_2(T\sm \{c\},-a_{n-2},-a_2))\in S_n(A)$. As
$$c+\sigma( U_1(T\sm \{c\},-a_{n-1}))\ \t{and}\ c+\sigma( U_2(T\sm \{c\},-a_{n-1},-a_{n-2}))$$
are on the list \eqref{list'} of all the elements of $S_n(A)$. Therefore, for some $x\in \{a_1,\ldots,a_{n-2}\}$ we have
$$c+\sigma( U_2(T\sm \{c\},-a_{n-2},-a_2))=c+\sigma( U_2(T\sm \{c\},-a_{n-1},-x)),$$
that is,
\begin{align*}
&d-(b+d)+\sum_{1<j<n-3}(b+jd)-(b+(n-3)d)+b+(n-2)d+c\\
=&\ d+\sum_{j=1}^{n-3}(b+jd)-b-(n-2)d+c-2x.
\end{align*}
It follows that $x=b$.
Combining this with the fact that $$x\in\{a_1,\ldots,a_{n-2}\}=\{d\}\cup\{b+jd:\ 0<j<n-2\},$$ we find that $x=d$.
Therefore $b=x=d$.

In view of the above, we have completed our proof of Theorem 1.4(i). \qed

\section{Proofs of Theorem 1.4(ii) and Theorem 1.5}
\setcounter{lemma}{0}
\setcounter{theorem}{0}
\setcounter{corollary}{0}
\setcounter{remark}{0}
\setcounter{equation}{0}

\begin{lemma}\label{lemma3} Let $n\gs 3$ be an integer.

{\rm (i)} Suppose that $A$ is a set of $2n$ integers with
\begin{equation}\label{A2n} A=\{c_1,c_2,\pm d\}\cup \{\pm (b+jd):\ j=1,\ldots n-2\},
\end{equation}
where $b,c_2,d\in\Z^+$,  $c_1\in \Z$, $c_1<b+(n-2)d $ and $c_2\gs b+(n-2)d$.
Then, for some $w\in S_n(A)\sm S_n(A\sm\{c_2\})$ we have
\begin{equation}\label{w_1} \min S_n(A\sm\{c_2\})<w<\max S_{n}(A\sm\{c_2\})+c_2-b-(n-2)d.
\end{equation}

{\rm (ii)} Suppose that $A$ is a set of $2n+1$ integers with
\begin{equation}\label{A2n'} A=\{c,\pm d\}\cup\{\pm (b+jd):\ j=1,\ldots ,n-1\},
\end{equation}
where  $b,c,d\in\Z^+$ and $c>b+(n-1)d$. Then, for some $w\in S_n(A)\sm S_n(A\sm\{c_2\})$ we have
\begin{equation}\label{w_2} \min S_n(A\sm\{c\})<w<\max S_{n}(A\sm\{c\})+c-b-(n-1)d.
\end{equation}
\end{lemma}
\Proof. Let $a_0=c_1$, $a_1=d,$ and $a_i=b+(i-1)d$ for $i=2,\ldots,n$.
\medskip

(i) We first prove part (i) of Lemma \ref{lemma3}.
Let $T=\{a_0,a_1,a_2,\ldots a_{n-1}\}$. Then
\begin{align*}
&a_0+\sigma( U_1(T\sm\{a_0\},-a_1))=a_0-a_1+\sum_{1<i<n}a_i
<\sum_{i=0}^{n-1} a_i<\sum_{i=0}^{n-2} a_i+c_2,
\end{align*}
and $a_0+\sigma( U_1(T\sm\{a_0\},-a_1))$ is the next term of $\sum_{i=0}^{n-1} a_i$
on the list \eqref{list'}.
Note also that
\begin{align*}
a_0+\sigma( U_1(T\sm\{a_0\},-a_1))<a_0-a_1+\sum_{1<i<n-1}a_i+c_2<\sum_{i=0}^{n-2} a_i+c_2.
\end{align*}

Now assume that there is no $w\in S_n(A)\sm S_n(A\sm\{c_2\})$ satisfying \eqref{w_1}. By the above, we must have
$$a_0-a_1+\sum_{1<i<n-1}a_i+c_2\in S_n(A\sm\{c_2\})$$
and hence
$$\sum_{i=0}^{n-1} a_i=a_0-a_1+a_2+\cdots +a_{n-2}+c_2.$$
It follows that
$$c_2=a_{n-1}+2a_1=b+nd.$$
Similarly, since
\begin{align*}
&a_0+\sigma( U_1(T\sm\{a_0\},-a_2))=a_0+a_1-a_2+a_3+\cdots +a_{n-1}\\
<&\ a_0+a_1-a_2+\cdots +a_{n-2}+c_2
\\<&\ a_0-a_1+\cdots +a_{n-2}+c_2=\sigma(T)
\end{align*}
and
\begin{align*}
&a_0+\sigma(U_1(T\sm\{a_0\},-a_2))=a_0+a_1-a_2+a_3+\cdots +a_{n-1}\\
<&\ \sigma(U_1(T\sm\{a_0\},-a_1))=a_0-a_1+a_2+\cdots +a_{n-1}<\sigma(T),
\end{align*}
we must have
$$a_0+a_1-a_2+\cdots +a_{n-2}+c_2=a_0-a_1+a_2+\cdots +a_{n-1}$$
and hence $b=d.$

Set
$$
\frac{A}{d}=\l\{\f ad:\ a\in A\r\}=\l\{\frac {c_1}d, \pm 1,\ldots , \pm (n-1), n+1\r\}.
$$
Then, all the elements of $S_n((A\sm\{c_2\})/d)$ have the same parity with $c_1/d +\sum_{i=1}^{n-1} i$.
Note that
$$w=\frac {c_1}d+1-2+3+\cdots+ (n-3)+n-1+n+1=\frac {c_1}d+\sum_{i=0}^{n-1} i-1$$
belongs to $S_n(A)\sm S_n(A\sm\{c_2\})$ and  the equality \eqref{w_1} holds.
So we get a contradiction.

(ii) Next we prove part (ii) of Lemma \ref{lemma3}.
In view of \eqref{A2n'}, we set
$$T=\{d\}\cup \{b+jd:\ j=1,\ldots n-1\}.$$
 Note that $u<\sigma (U_1(T,-a_2))<\sigma(U_1(T,-a_1))<\sigma (T)$ for all $u\in S_n(A\sm\{c\})\sm\{\sigma (U_1(T,-a_2)),\sigma (U_1(T,-a_1)),\sigma (T)\}$. Also,
$$\sigma(U_1(T,-a_1))=-a_1+\sum_{1<i\ls n}a_i<\sigma(T)= \sum_{i=1}^{n} a_i<\sum_{i=1}^{n-1} a_i+c$$
and
$$\sigma(U_1(T,-a_1))<-a_1+\sum_{1<i< n}a_i+c<\sum_{i=1}^{n-1} a_i+c.$$

Now assume that there is no $w\in S_n(A)\sm S_n(A\sm\{c\})$ satisfying \eqref{w_2}. By the last paragraph, we must have
$$\sum_{i=1}^{n} a_i=-a_1+\sum_{1<i<n}a_i+c.$$
It follows that
\begin{equation}\label{cbd}
c=2a_1+a_n=b+(n+1)d.
\end{equation}
Similarly, since
\begin{align*}
&\sigma(U_1(T,-a_2))=a_1-a_2+a_3+\cdots +a_n\\
<&\ a_1-a_2+\cdots +a_{n-1}+c\\
<&\ \sigma(T)=-a_1+\cdots +a_{n-1}+c.
\end{align*}
and
$$
\sigma( U_1(T,-a_2))<\sigma( U_1(T,-a_1))<\sigma(T),
$$
we must have
$$a_1-a_2+\cdots a_{n-1}+c=-a_1+a_2+\cdots +a_{n},$$
and hence
$$c=-2a_1+2a_2+a_n.$$
Combining this with \eqref{cbd}, we get $b=d.$
Now we see that
$$
\frac{A}{d}=\{\pm 1,\ldots , \pm n, n+2\}.
$$
As all the elements of $S_n((A\sm\{c\})/d)$ have the same parity with $\sum_{i=1}^n i$, and
$$w=1-2+3+\cdots +n-2+n+n+2=\sum_{i=1}^n i-1$$
belongs to $S_n(A)\sm S_n(A\sm\{c\})$, we get a contradiction in view of the equality \eqref{w_2}.

By the above, we have completed the proof of Lemma \ref{lemma3}.
\qed

\medskip
\noindent{\it Proof of Theorem} 1.4(ii). The ``if" direction is easy since we can verify the result directly.

Now assume the equality \eqref{A=} with $|A|=k=2n$.
By Lemma 3.1, $A$ contains at least $n-1$ negative integers and also at least $n-1$ positive integers.

{\tt Step I}. Prove that
\begin{equation}\label{+-} |\{a\in A:\ a>0\}|=n=|\{a\in A:\ a<0\}|.
\end{equation}

Suppose that \eqref{+-} fails. Then either $A$ or $-A$ contains exactly $n+1$ nonnegative integers.
Since $|S_n(-A)|=|S_n(A)|$, without loss of generality, we may assume that $|A\cap \N|=n+1$.
Write all the $2n$ elements of $A$ in the ascending order
$$b_{n-1}<b_{r-1}<\cdots <b_1<a_1<\cdots <a_{n+1}$$
with $b_1<0\ls a_1$. Note that
$$|S_n(\{a_1,\ldots,a_{n+1}\})|=|n^\wedge\{a_1,\ldots,a_{n+1}\}|=n+1.$$
Also, $a_1+\ldots+a_n$ is the only common element of $S_n(A\sm\{a_{n+1}\})$ and $S_n(\{a_1,\ldots,a_{n+1}\})$.
Thus
$$|S_n(A\sm\{a_{n+1}\})|=\frac {n(n+1)}2+1-n=\frac {n(n-1)}2+1$$
in the spirit of the inequality \eqref{Sn}.

By Theorem 1.4(i), we have
$$A\sm\{a_{n+1}\}=\{c,\pm d\}\cup\{\pm (b+jd):\ j=1,\ldots n-2\}$$ for some $b,d\in \Z^+$ and $c\in \Z$.
\medskip

{\it Case} 1. $c<b+(n-2)d$.
 \medskip

 By Lemma \ref{lemma3}(i), for certain $w\in S_n(A)\sm S_n(A\sm\{a_{n+1}\})$ we have $w<a_1+\cdots +a_{n+1}-b-(n-2)d$. Note that $w\notin S_n(\{a_1,\ldots,a_{n+1}\})$ since
 $$a_1+\cdots +a_{n+1}-b-(n-2)d=\min (S_n(\{a_1,\ldots,a_{n+1}\})\sm \{a_1+\ldots+a_n\}).$$
Hence
\begin{align*}
|S_n(A)|\gs&\ |S_n(A\sm\{a_{n+1}\})|+|S_n(\{a_1,\ldots,a_{n+1}\})|+|\{w\}|-1
\\=&\ \frac {n(n-1)}2+2,
\end{align*}
which contradicts the equality \eqref{A=}.

\medskip

{\it Case} 2. $c>b+(n-2)d$.
 \medskip

Clearly, $a_1=d$, $a_i=b+(i-1)d$ for all $i=2,\ldots,n-1$, and $a_{n+1}>a_n>b+(n-2)d=a_{n-1}$.
Let $T=\{a_1,a_2,\ldots ,a_n\}$.
Then the $n(n+1)/2+1$ numbers
\begin{equation}\label{list-}\begin{aligned}&\sum_{i=2}^{n+1} a_i,\ a_1+\sum_{i=3}^{n+1} a_i,\ \cdots ,\ \sum_{i=1}^{n} a_i=\sigma (T),
\\&a_n+\sigma(U_1(T\sm \{a_n\},-a_i))\ (i=1,\ldots,n-1),
\\&a_n+\sigma(U_2(T\sm \{a_n\},-a_{n-1},-a_j))\ (j=1,\ldots,n-2),
\\&\cdots,\ a_n+\sigma(U_{n-1}(T\sm \{a_n\},\ -a_{n-1},-a_{n-2}\ldots ,-a_1))
\end{aligned}
\end{equation}
are distinct elements of $S_n(A)$ and they are in descending order.
Note that
\begin{align*}
a_1+\sum_{i=3}^{n+1} a_i >a_1+a_2+\sum_{i=4}^{n+1} a_i >\cdots>a_n+\sigma(U_1(T\sm \{a_n\},-a_1)).
\end{align*}
are $n+1$ consecutive terms on the list \eqref{list-} in descending order.
On the other hand, the $n+1$ numbers
\begin{align*}
&a_1+\sum_{i=3}^{n+1} a_i >-a_1+\sum_{i=3}^{n+1} a_i \\
>&-a_1+a_2+\sum_{i=4}^{n+1} a_i>\cdots >a_n+\sigma(U_1(T\sm \{a_n\},-a_1)).
\end{align*}
in $S_n(A)$ are also in descending order. Therefore
$$a_1+a_2+\sum_{i=4}^{n+1} a_i=-a_1+\sum_{i=3}^{n+1} a_i$$
and hence $a_1+a_2=-a_1+a_3$. But
$$2a_1+a_2-a_3=2d+b+d-(b+2d)=d\not=0,$$
so we get a contradiction.
\medskip

{\tt Step II}. Show that $A$ has the desired form.

List all the $2n$ elements of $A$ in the ascending order:
$$b_{n}<\cdots <b_1<a_1<\cdots <a_n,$$
where $b_1<0<a_1$.  Set $T=\{a_1,\ldots ,a_n\}.$
Recall that the list \eqref{list} gives $r(r+1)/2+1$ elements of $S_n(A)$ in descending order.
By similar arguments, the $n(n+1)/2+1$
\begin{equation}\label{list2n}\begin{aligned}&\sigma(T),\sigma(U_1(T,b_i))\ (i=1,\ldots,n),
\\&\sigma(U_2(T,b_n,b_i))\ (i=1,\ldots,n-1),
\\&\cdots,\sigma(U_n(T,b_n,b_{n-1}\ldots ,b_1))
\end{aligned}
\end{equation}
are distinct elements of $S_n(A)$ and they are in descending order.

We now claim that $-b_j\in A$ for all $1\ls j\ls n$.
Assume on the contrary that $-b_j\notin A$ for certain $j\in\{1,\ldots,n\}$.
Then
$$|S_n(\{b_j,a_1,\ldots,a_{n}\})|=|n^\wedge\{b_j,a_1,\ldots,a_{n}\}|=n+1.$$
Also, $b_j+\sum_{i=1}^{n-1} a_i$ is the only common element of $S_n(A\sm\{a_{n}\})$ and $S_n(\{b_j,a_1,\ldots,a_{n}\})$.
Thus
$$
|S_n(A\sm\{a_n\})|=\frac {n(n+1)}2+1-n-1+1=\f{n(n-1)}2+1
$$
in the spirit of the inequality \eqref{Sn}.
By Theorem 1.4(i), there are $b,d\in\Z^+$ such that
$A\sm\{a_n\}=\{b_j,\pm d\}\cup \{\pm (b+jd):\ j=1,\ldots n-2\}.$
Hence
$$A=\{b_j,a_n,\pm d\}\cup \{\pm (b+jd):\ j=1,\ldots n-2\}.$$
By Lemma \ref{lemma3}(i), for some $w\in S_n(A)\sm S_n(A\sm\{a_n\})$ we have $w<b_j+a_n+d+\sum_{i=1}^{n-3} (b+id)$,
which leads to a contradiction.

By the claim in the last paragraph, we have $a_i=-b_i$ for all $i=1,\ldots,n$.

Recall that all the elements of $S_n(A)$ are exhibited on the list \eqref{list2n}, and
\begin{equation}\label{2n1}\begin{aligned}
\sigma(U_1(T,b_i))\ (i=1,\ldots,n),\ \sigma(U_2(T,b_n,b_1))
\end{aligned}
\end{equation}
are $n+1$ consecutive terms on the list \eqref{list2n} in descending order. Note also that the $n+1$ numbers
\begin{equation}\label{2n2}\begin{aligned}
\sigma(U_1(T,b_i))\ (i=1,2),\ \sigma(U_2(T,b_j,b_1))\ (j=2,\ldots ,n)
\end{aligned}
\end{equation}
in descending order also belong to $S_n(A).$
As \eqref{2n1} and \eqref{2n2} have the same first term and the same last term,
the two lists \eqref{2n1} and \eqref{2n2} are identical.
Therefore, for any $2\ls i\ls n$, we have
$$\sigma(U_1(T,b_{i+1}))=\sigma( U_2(T,b_{i},b_1))$$
and hence $a_{i+1}-a_i=a_1$.
Note that $d=a_1$ and $b=a_2-a_1$ are positive integers. For each $i=2,\ldots,n-1$ we have $-b_i=a_i=b+(i-1)d$.
Thus
$$A=\{\pm d,\pm (b+jd):\ j=0,1,\ldots n-1\}.$$

Now assume  $n\gs4$. We want to prove that $b=d$. Note that
\begin{align*}
&\ \sigma( U_2(T,-a_{n-1},-a_2))\\
=&\ d-(b+d)+\sum_{1<i<n-2}(b+id)-(b+(n-2)d)+b+(n-1)d,
\end{align*}
is smaller than $\sigma(U_1(T,-a_n))$ and greater than
$\sigma( U_2(T,-a_n,-a_{n-1}))$.
As $\sigma( U_1(T,-a_n))$ and $\sigma( U_2(T,-a_n,-a_{n-2}))$
are on the list \eqref{list2n} of all the elements of $S_n(A)$, for some $x\in \{a_1,\ldots,a_{n-1}\}$ we have
$$\sigma( U_2(T,-a_{n-1},-a_2))=\sigma( U_2(T,-a_n,-x)),$$
that is,
\begin{align*}
&\ d-(b+d)+\sum_{1<i<n-2}(b+id)-(b+(n-2)d)+b+(n-1)d\\
=&\ d+\sum_{1\ls i\ls n-2}(b+id)-(b+(n-1)d)-2x
\end{align*}
It follows that $x=b$.
Combining this with the fact that $$x\in\{a_1,\ldots,a_{n-1}\}=\{d\}\cup\{b+jd:\ 0<j<n-1\},$$ we find that $x=d$.
Therefore $b=x=d$.

In view of the above, we have completed our proof of Theorem 1.4(ii). \qed
\medskip

\noindent{\bf Proof of Theorem 1.5.}
By Lemma 3.1, $A$ contains at least $n-1$ negative integers and at least $n-1$ positive integers.
Since $|S_n(A)|=|S_n(-A)|$, without loss of generality, we may assume
$$|\{a\in A:\ a\gs0\}|\gs |\{a\in A:\ a<0\}|.$$
 Thus, either $A$ contains exactly $n-1$ negative numbers, or $A$ contains at least $n$ negative numbers and $n+1$ nonnegative numbers.
\medskip

{\it Case} 1. $A$ contains exactly $n-1$ negative numbers.
\medskip

Let $A'$ be a set consisting of $2n$ consecutive terms of $A$ such that it contains all negative numbers in $A$.
In view of Theorem 1.2,
$$|S_n(A')|\gs \frac {n(n+1)}2 +1.$$
 As $A'\neq \{\pm d\}\cup\{\pm (b+jd):\ j=0,\ldots,n-1\}$ for any $b,d\in\Z^+$, by Theorem 1.4(ii) we have $|S_n(A')|>n(n+1)/2 +1.$

 Let us list all the elements of $A'$ in the ascending order:
$$b_{n-1}<\cdots <b_1<a_1\cdots <a_{n+1}.$$
Set $B= \{a_2,a_2,\cdots a_{n+1} \}.$
Clearly, both  $S_n(A')$ and $n^\wedge(B\cup (A\sm A'))$ are subsets of $S_n(A)$.
Note that $a_2+\cdots + a_{n+1}$ is the only common element of $S_n(A')$ and $n^\wedge(B\cup (A/A'))$.
Thus
\begin{align*}
S_n(A) &\gs |S_n(A')|-1+|n^\wedge(B\cup (A/A'))|\\
&\gs\frac {n(n+1)}2+(k-n)n-n^2+1+1\\
&= kn-\frac 32 n^2+\frac 12 n+2.
\end{align*}

\medskip
{\it Case} 2. $A$ contains at least $n$ negative numbers and at least $n+1$ nonnegative numbers .
\medskip

Let $A'$ be a set consisting of $2n$ consecutive terms of $A$
such that
$$|\{a\in A':\ a<0\}|=n=|\{a\in A':\ a\gs0\}|.$$
Set $B=A\cap \N$ and $C=\{a\in A:\ a<0\}$ with $|B|=t>n$ and $|C|=k-t\gs n$.

By (1.1) we have
$$|n^\wedge B|=|S_n(B)| \gs nt-n^2+1\ \t{and}\ |n^\wedge C|=|S_n(C)| \gs n(k-t)-n^2+1.$$
Note that
$$ |n^\wedge B\cap S_n(A)|=1,\  |n^\wedge C\cap S_n(A)|=1,\
\t{and}\  n^\wedge B\cap n^\wedge C=\emptyset.$$
In the spirit of (1.2), we also have
$$|S_n(A')|\gs (2n-1)n-3\bi n2+1=\f{n(n+1)}2+1.$$
Thus
\begin{align*}
|S_n(A)|\gs &\ |n^\wedge B|+|S_n (A')|+|n^\wedge C|-2\\
\gs &\ nt-n^2+1+\l(\frac {n(n+1)}2+1\r)+n(k-t)-n^2+1-2\\
=&\ kn-\frac 32 n^2+\frac 12 n+1.
\end{align*}

Suppose that $|S_n(A)|=kn-\frac 32 n^2+\frac 12 n+1.$ Then
$$|n^\wedge B|=nt-n^2+1,\ |S_n (A')|=\frac {n(n+1)}2+1,\ |n^\wedge C|=n(k-t)-n^2+1.$$
Applying Theorem 1.4(ii), we see that
$$A'=\{\pm d,\pm (b+jd):\ j=1,\ldots n-1\}$$
for some $b,d\in\Z^+$. Let $B'=(A\sm A')\cap \N$.
By Lemma 4.1(ii), there is a number $w\in S_n(A'\cup \{\min B'\})\sm S_n(A')$ such that $\min S_n(A')<w<\max S_n(A')+\min B'-b-(n-1)d,$ and it is easy to see $w\notin n^\wedge B$. Hence
\begin{align*}
|S_n(A)|\gs &|n^\wedge B|+|S_n (A')|+|n^\wedge C|+1-2
\\\ =&kn-\frac 32 n^2+\frac 12 n+2.
\end{align*}
So we get a contradiction.

In view of the above, we have completed our proof of Theorem 1.5.
\qed

\setcounter{conjecture}{0}
\end{document}